\documentclass[10pt]{amsart}

\usepackage{a4wide,bm}
\usepackage{geometry}
\usepackage{thmtools}
\usepackage{thm-restate}
\usepackage[colorlinks=true]{hyperref}
\hypersetup{colorlinks=true,  linkcolor=blue,  filecolor=blue,  citecolor=blue, urlcolor=blue}
\usepackage[capitalize,noabbrev]{cleveref}
\usepackage{ wasysym }
\usepackage{amssymb}
\usepackage{euscript}
\usepackage{xcolor}
\usepackage{indentfirst}
\usepackage{microtype}
\usepackage{graphicx}
\usepackage{mathtools}
\usepackage{amsmath}
\usepackage{subcaption}
\usepackage{float}
\usepackage{enumitem}  
\usepackage{mathptmx}

\usepackage{mathrsfs}
\usepackage[all]{xy}
\usepackage[ruled,vlined]{algorithm2e}
\usepackage[T1]{fontenc}
\usepackage{textcomp}
\usepackage{tikz}
\usepackage{capt-of}
\usepackage{lineno}
\usepackage{bm}
\usepackage{comment}
\usepackage{titlesec}
\usepackage{listings}

\theoremstyle{plain}
    \newtheorem{theorem}{Theorem}[section]
    \newtheorem{lemma}[theorem]{Lemma}
    \newtheorem{proposition}[theorem]{Proposition}
    
\theoremstyle{definition}
    \newtheorem{remark}[theorem]{Remark}

    \newtheorem{definition}[theorem]{Definition}
    
    \newtheorem{question}[theorem]{Question}

\DeclareMathOperator{\Gal}{Gal}

\DeclareMathOperator{\GL}{GL}

\newcommand{\R}{\mathbb{R}}
\newcommand{\Q}{\mathbb{Q}}
\newcommand{\Z}{\mathbb{Z}}

\renewcommand{\P}{\mathbb{P}}
\newcommand{\C}{\mathbb{C}}

\newcommand{\F}{\mathbb{F}}

\newcommand{\SF}{\mathcal{S}(F)}

\newcommand*\Zp{\Z_p}
\newcommand*\Qp{\Q_p}

\newcommand*\Qpb{\overline{\Q}_p}

\newcommand{\norm}[1]{\left\lVert#1\right\rVert}

\newcommand*\varhrulefill[1][0.4pt]{\leavevmode\leaders\hrule height#1\hfill\kern0pt}
\newcommand*{\Scale}[2][4]{\scalebox{#1}{$#2$}}%

\makeatletter
\renewcommand{\section}{\@startsection {section}{1}{\z@}%
	{-3.5ex \@plus -1ex \@minus -.2ex}%
	{2.3ex \@plus .2ex}%
	{\large \scshape \bfseries \filcenter}}
\renewcommand{\@secnumfont}{\bfseries}
\makeatother

\title{Lines on $p$-adic and real cubic Surfaces}

\author[R.~Ait El Manssour]{{Rida Ait El Manssour}}
\address{Rida Ait El Manssour (MPI MiS)}
\email{rida.manssour@mis.mpg.de}

\author[Y.~El Maazouz]{{Yassine El Maazouz}}
\address{Yassine El Maazouz (UC Berkeley)}
\email{yassine.el-maazouz@berkeley.edu}

\author[E.~Kaya]{{Enis Kaya}}
\address{Enis Kaya, (KU Leuven)}
\email{enis.kaya@kuleuven.be}

\author[K.~Rose]{{Kemal Rose}}
\address{Kemal Rose (MPI MiS)}
\email{kemal.rose@mis.mpg.de}


\begin{document}
	
	\begin{abstract}  
		We study lines on smooth cubic surfaces over the field of $p$-adic numbers, from a theoretical and computational point of view. Segre showed that the possible counts of such lines are $0,1,2,3,5,7,9,15$ or $27$. We show that each of these counts is achieved. Probabilistic aspects are investigated by sampling both $p$-adic and real cubic surfaces from different distributions and estimating the probability of each count. We link this to recent results on probabilistic enumerative geometry. Some experimental results on the Galois groups attached to $p$-adic cubic surfaces are also discussed.
	\end{abstract}
	
\maketitle
\setcounter{tocdepth}{2}
%\tableofcontents

\section{Introduction}\label{Sec:1}

Smooth cubic surfaces have been intensively studied by geometers, and some of the basic sources for these objects are \cite{Hartshorne1977,Manin1986,Dolgachev2012}. The problem of counting lines on such surfaces also has a long standing history and there is a vast literature on the topic. In the nineteenth century, Cayley \cite{Cayley1849} and Schl{\"a}fli \cite{Schlafli1858} proved that a smooth cubic surface over the complex numbers always has $27$ lines, while over the real numbers it can have $3,7,15$ or $27$ lines. Segre \cite{Segre1949} showed that over any field, a smooth cubic surface can only have $0,1,2,3,5,7,9,15$, or $27$ lines\footnote{In his work, Segre actually points out that his statement fails in characteristic $2$, but this is not correct; see \cite[Section~3.1]{McKean}.}. Depending on the base field, some of these numbers may not be attained. Since then, the problem of counting lines on cubic surfaces over different fields attracted many mathematicians and some cases have already been classified; see Table~\ref{PossibleLineCounts}.

\begin{table}[H] 
    \centering
        \begin{tabular}{|c|c|c|}
                \hline
                Base field & Numbers attained & Reference \\[0.6ex]
                \hline
                $\C$ & 27 & \cite{Cayley1849} \\[0.6ex]
                \hline
                $\R$  & 3, 7, 15, 27 & \cite{Schlafli1858} \\[0.6ex]
                \hline
                $\Q$ & 0, 1, 2, 3, 5, 7, 9, 15, 27 & \cite{Segre1949} \\[0.6ex]
                \hline
                $\F_{q}$, $q > 5$ \text{odd} & 0, 1, 2, 3, 5, 7, 9, 15, 27 &  \cite{LoughranTrepalin} \\[0.6ex]
                \hline
                $\F_{2^n}$, $n\geq 2$  & 0, 1, 2, 3, 5, 7, 9, 15, 27 & \cite{LoughranTrepalin} \\[0.6ex]
                \hline
                $\F_{5}$ & 0, 1, 2, 3, 5, 7, 9, 15 & \cite{LoughranTrepalin}\\[0.6ex]
                \hline
                $\F_2, \ \F_3$ & 0, 1, 2, 3, 5, 9, 15 & \cite{Dickson1914/15, LoughranTrepalin}\\[0.6ex]
                \hline
        \end{tabular}
        \caption{Classification of possible line counts for the fields $\C$, $\R$, $\Q$, $\F_q$.}
        \label{PossibleLineCounts}
\end{table}

Recently, McKean \cite{McKean} showed that all of the numbers $0,1,2,3,5,7,9,15$ and $27$ occur when the base field is finitely generated with at least $22$ elements or a finite transcendental extension of an arbitrary field; see \cite[Corollary~1.6]{McKean}. 

In this text, our main aim is to get an idea of how many lines are on a cubic surface over $\Q_p$ or $\R$ from a probabilistic point of view. To this end, we conduct numerical experiments when the base field is $\Q_7$ or $\R$. While over the real numbers, a cubic surface can only have either $3,7,15$ or $27$ lines, by adapting the approach in an earlier version of \cite{McKean}\footnote{This was done prior to the appearance of McKean's stronger result \cite[Theorem 1.3]{McKean}.} to the case we are interested in, we get the following:

\begin{theorem}\label{thm:1}
    Let $n\in\{0, 1, 2, 3, 5, 7, 9, 15, 27\}$. Then there exists a smooth cubic surface over $\Qp$ that contains exactly $n$ lines. In other words, all possible line counts occur when the base field is $\Qp$.
\end{theorem}

In \cref{Sec:2}, we explain how to explicitly construct (by blowing up the projective plane in $6$ suitable points) a smooth cubic surface having any of the line counts mentioned in \cref{thm:1}. \cref{Sec:3} and \cref{Sec:4} focus on probabilistic computations and heuristics for both the $p$-adic and real case, respectively. For the $p$-adic case, we sample from the family of smooth cubic surfaces with four different probability measures (see \cref{Sec:3}), and compute the probability of seeing each number of lines. \cref{tab:p-adicResults} summarizes the distributions of the number of lines when $p = 7$. For the real case, we consider a one-parameter family $(\P_\lambda)_{0< \lambda < 1}$ of Gaussian distributions studied in \cite{Kost90}. The probability distribution of line counts is then a curve in the $3$-simplex which is depicted in Figure~\ref{fig:curveReals}. The Galois groups attached to smooth cubic surfaces are of special interest. The final section of this paper (\cref{Sec:5}) is devoted to experimental results on which Galois groups appear for cubic surfaces defined over $\Qp$. Our results suggest that Galois groups should be quite small and usually abelian. This motivates an inverse Galois problem for smooth cubic surfaces over $\Q_p$ in line with Elsenhans and Jahnel's work \cite{ElsenhansJahnel2015} on the inverse Galois problem for smooth cubic surfaces over $\Q$.

\medskip

Most of the results in this text were found by computation. The codes and data are made available at
\begin{equation}\label{eq:link_for_code}
    \text{\url{https://mathrepo.mis.mpg.de/27pAdicLines/index.html}}.
\end{equation}
Our computations were carried out using the computer algebra systems Magma \cite{Magma}, Julia \cite{Julia} and Macaulay2 \cite{Macaulay2}.

\vspace{1cm}

\begin{center}
    \textbf{Acknowledgements}
\end{center}
We thank Claus Fieker, Tom Fisher, Stevan Gajovi\'c, Marta Panizzut, Emre Sert\"oz and Bernd Sturmfels for valuable discussions. We also thank Avinash Kulkarni, Daniel Loughran and Jean Pierre Serre for their valuable comments on earlier versions of this manuscript. Finally, we thank the referee for
suggestions that improved the paper.

\section{Lines on Smooth \texorpdfstring{$p$-}~adic Cubic Surfaces} \label{Sec:2}

In the following, we recall a well-known representation of smooth cubic surfaces: the blow-up of the projective plane $\P^2$ at six $\Qpb$-rational points in general position is a smooth cubic surface (see, for example, \cite{OctanomialModel} for a detailed treatment). Up to a change of coordinates, the six points lie on the cuspidal cubic curve
\[ 
    \mathcal{C}:  \theta \longmapsto (1:\theta:\theta^3)
\]
and can be represented as the columns of the matrix
\begin{equation}
\label{eq: cuspidal}
\begin{bmatrix}
1 & 1 & 1 & 1 & 1 & 1\\
\theta_1 & \theta_2 & \theta_3 & \theta_4 & \theta_5 & \theta_6\\
\theta_1^3 & \theta_2^3 & \theta_3^3 & \theta_4^3 & \theta_5^3 & \theta_6^3\\
\end{bmatrix}.
\end{equation}
The six points are said to be in \emph{general position} if no three of them lie on a line and not all of them lie on a conic. Equivalently the points are in general position if both the maximal minors of \eqref{eq: cuspidal}, which are given by
\[
    \lbrack i j k \rbrack  \coloneqq (\theta_i - \theta_j )(\theta_i - \theta_k)(\theta_j - \theta_k)(\theta_i + \theta_j + \theta_k) \quad \text{ for }  1 \leq i < j < k \leq 6,
\]
and the polynomial
\[
  \Scale[0.95]{\lbrack 1 3 4 \rbrack \lbrack 1 5 6 \rbrack \lbrack 235 \rbrack \lbrack 246 \rbrack  - \lbrack 135 \rbrack \lbrack 146 \rbrack \lbrack 234 \rbrack \lbrack 256 \rbrack 
\coloneqq  (\theta_1 + \theta_2 + \theta_3 + \theta_4 + \theta_5 + \theta_6) \prod_{1 \leq i < j \leq 6} (\theta_i - \theta_j)}  ,
\]
 do not vanish. ~These polynomials split into linear forms that together determine the hyperplane arrangement
of type $E_6$ (see \cite[Section~6]{RenQingchun2014} for more details):
    \begin{align} \label{eq: arrangement}
        \begin{split}
            \theta_i  -  \theta_j \hspace{28pt}
            & \text{ for } 1 \leq i < j \leq 6, \\
            \theta_i + \theta_j + \theta_k \hspace{15pt}
            & \text{ for } 1 \leq i < j < k \leq 6, \\
            \theta_1 + \theta_2 + \cdots + \theta_6.
        \end{split}
    \end{align}
\begin{definition}\label{def: octanomial surface}
    Let $F = \prod\limits_{i = 1}^6 (X - \theta_i) \in \Qp[X]$ be a univariate polynomial of degree $6$ whose roots lie in the complement of the hyperplane arrangement \eqref{eq: arrangement}.
We denote by $\SF$ the smooth cubic surface that is the blow-up of the projective plane at the six points $\{[1: \theta_i: \theta_i^3]\}_{1 \leq i \leq 6}$.

The defining equation of $\SF$ in $\P^3$ in terms of the $\theta_i$ is determined in \cite[Equation~4]{OctanomialModel}. The $27$ $\Qpb$-rational lines on $\SF$ are of three distinct types:
\begin{enumerate}[label=(\roman*)]
    \item $\{E_i \colon i =1,\dots,6\}$, where $E_i$ is the exceptional divisor of the point $\mathcal C (\theta_i)$;
    \item $\{F_{i, j} \colon i,j =1,\dots,6,\ i\neq j\}$, where $F_{i, j}$ is the strict transform of the line passing through the points $\mathcal C (\theta_i)$ and $\mathcal C (\theta_j)$; 
    \item $\{G_i \colon i =1,\dots,6\}$, where $G_i$ is the strict transform of the unique conic passing through the points $\{\mathcal C (\theta_1),\dots,\mathcal C (\theta_6)\} \setminus \{\mathcal C (\theta_i)\}$.
\end{enumerate}

\end{definition}
 
In view of constructing smooth cubic surfaces that have a prescribed number of $\Qp$-rational lines, we investigate the action of the absolute Galois group $G = \Gal(\Qpb/\Qp) $ on $\SF$, and show that the line-count depends only on the decomposition of $F$ into irreducible factors over $\Qp$.
\begin{lemma}\label{lem:NumberOfLines}
    The smooth cubic surface $\SF$ is $\Qp$-rational and contains exactly
    \[
        2\ell + q + \binom{\ell}{2}
    \]
    $\Qp$-rational lines, where $\ell$ and $q$ are the number of linear and quadratic irreducible factors of $F$ in $\Qp[X]$, respectively.
\end{lemma}
\begin{proof}
The absolute Galois group $G$ acts on $\SF$ and lines on it by permuting the roots $\theta_1, \dots, \theta_6$. In particular, for every element $\sigma$ of $G$, we have
\[
\sigma(\SF) = \mathcal{S}(\sigma(F)), \ \sigma(E_i) = E_{\sigma(i)}, \ \sigma(F_{i, j}) = F_{\sigma(i), \sigma(j)}, \ \sigma(G_i) = G_{\sigma(i)},
\]
where we define $\sigma(i)$ by enforcing $\sigma(\theta_i) = \theta_{\sigma(i)}$. Since $G$ acts trivially on $F \in \Qp[X]$, the surface $\SF$ is stable under $G$. Hence it is $\Qp$-rational. Note that the number of $G$-stable $E_i$'s is $\ell$, the number of $G$-stable $F_{i,j}$'s is $q + \binom{\ell}{2}$, and the number of $G$-stable $G_i$'s is $\ell$. This finishes the proof.
\end{proof}

\begin{proof}[Proof of Theorem~\ref{thm:1}]
This is a direct consequence of the above lemma. The polynomials in Table~\ref{9polynomials} have the desired numbers of irreducible linear and quadratic factors. The discriminant, the degree $5$ coefficients $\theta_1 + \dots +\theta_6$ and the expression
\[\prod_{ \# \{i, j, k\} = 3    } \theta_i + \theta_j + \theta_k\]
are symmetric and can be expressed as polynomials in the coefficients of $F$. The non-vanishing of these quantities for the polynomials in \cref{9polynomials} was certified using Macaulay2 (see \eqref{eq:link_for_code}); therefore, the equations in $\eqref{eq: arrangement}$ are not satisfied for any of these $9$ polynomials.
\begin{table}[H]
    \centering
            \begin{tabular}{|c|c|c|}
                \hline
                Polynomial $F$ & $(\ell,q)$& $\#$ of lines on $\SF$ \\[0.6ex]
                \hline
                $X^6 + pX^5 + p$ & $(0, 0)$ & $0$ \\[0.6ex]
                \hline
                $(X^4 + p)(X^2 + pX + p)$ & $(0, 1)$ & $1$ \\[0.6ex]
                \hline
                $X(X^5 + pX^4 + p)$& $(1,0)$ & $2$ \\[0.6ex]
                \hline
                $(X^2 + p)(X^2 + pX + p)(X^2 + p^2X + p)$ & $(0, 3)$ & $3$ \\[0.6ex]
                \hline
                $(X + 1)(X + 2)(X^4 + p)$ & $(2, 0)$ & $5$ \\[0.6ex]
                \hline
                $(X + 1)(X + 2)(X^2 + p)(X^2 + pX + p)$ & $(2, 2)$ & $7$ \\[0.6ex]
                \hline
                $(X + 1)(X + 2)(X + 3)(X^3 + pX^2 + p)$ & $(3, 0)$ & $9$ \\[0.6ex]
                \hline
                $(X + 1)(X + 2)(X + 3)(X + 4)(X^2 + p)$ & $(4, 1)$ & $15$ \\[0.6ex]
                \hline
                $X(X + 1)(X + 2)(X + 3)(X + 4)(X + 5)$ & $(6, 0)$ & $27$ \\[0.6ex]
                \hline
            \end{tabular}
        \caption{A list of surfaces $\SF$ that realize all possible line-counts over$~\Qp$.}
        \label{9polynomials}
    \end{table}
\end{proof}

\begin{remark}\label{rem:DoubleSix}

    \begin{enumerate}
    
        \item It is clear that \cref{thm:1} continues to hold if we replace $\Q_p$ with any finite extension of $\Q_p$. In fact, it holds for any local field (one just needs to modify the above proof slightly). But, as we mentioned in the introduction, this is a special case of a theorem proved by McKean, after our work was completed, in the final version of his paper; see \cite[Corollary~1.6]{McKean}. We think it is still valuable to have explicit polynomials that yield the desired surfaces.
    
        \item It is important to note that, while any cubic surface over $\overline{\Q}_p$ is isomorphic to the blow-up of $6$ points in general position in $\overline{\Q}_p\P^2$ that is not the case over $\Q_p$. In other words, there are cubic surfaces defined over $\Q_p$ which do not arise from the blow-up construction as in \cref{Sec:2}. The cubic surfaces arising in that way are the cubic surfaces which have a Galois invariant (i.e. $\Gal(\overline{\Q}_p/ \Q_p)$-invariant) \emph{double-six}. A double-six is a configuration of $12$ lines $\{L_1, \dots, L_6\}$ and $\{L_1', \dots, L_6'\}$ in $\overline{\Q}_p\P^3$ arranged in matrix form 
        \[
            \begin{pmatrix}
                L_1  & L_2  & L_3  & L_4  & L_5  & L_6\\
                L_1' & L_2' & L_3' & L_4' & L_5' & L_6'
            \end{pmatrix},
        \]
        such that any two lines of these $12$ are secant if and only if they are not on the same row or column. With the notation of \cref{def: octanomial surface}, the configuration of lines $\{E_1, \dots, E_6\}$ and $\{G_1, \dots, G_6\}$ is a Galois invariant double-six; see, for example, \cite[Remark~V.4.9.1]{Hartshorne1977}.
        
    \end{enumerate}
\end{remark}

\section{Heuristics for the \texorpdfstring{$p$-}~adic Numbers} \label{Sec:3}

There are plenty of ways one can construct a smooth cubic surface over $\Q_p$. Each leads to a different way of sampling such surfaces and hence a measure on the space of smooth cubic surfaces. For our probabilistic investigations, we focus on four sampling methods and describe both the measures and how the sampling process is implemented in our computations.

\subsection{The Haar Measure}

    We endow the space of cubics $\mathbb{Q}_p[ x_0 , x_1, x_2, x_3]_{( 3 )}$ with its natural measure, defined by choosing the $20$ coefficients $\xi_{\alpha}$ of the degree $3$ cubic
    \begin{equation} \label{eq:cubicpoly}
            f = \sum_{|\alpha| = 3} \xi_\alpha x^{\alpha}    
    \end{equation}
    to be independent and uniformly distributed in $\Zp$ with respect to the Haar measure on $\Z_p$. 
    
    \begin{definition} \label{def:Haarmeasure}
        The probability measure of the random smooth cubic surface $\mathcal{S}(f)$ defined as the zero set of a random polynomial $f$ as in \eqref{eq:cubicpoly} is called the \emph{Haar measure} on the space of cubic surfaces.
    \end{definition}

    Note that the singular cubics lie on a hypersurface of $\mathbb{Q}_p[ x_0 , x_1, x_2, x_3]_{( 3 )}$ and hence with probability $1$ the random cubic $f$ in \eqref{eq:cubicpoly} defines a smooth surface. Moreover, the Haar measure in \cref{def:Haarmeasure} is invariant under change of variables, i.e.
    \[
    \text{for any } g \in \GL_4(\Z_p),\ f(g\cdot x)\text{ has the same distribution as } f(x).
    \]

    \begin{remark}
        We remark that one may make a different choice of basis in \eqref{eq:cubicpoly} and get another measure on the space of polynomials $\Q_p[x_0, x_1, x_2, x_3]_{(3)}$. The reason we chose the monomial basis in \eqref{eq:cubicpoly} is because it is guaranteed to be invariant under change of variables, moreover when $p$ is big enough, namely $p > 3$, the resulting measure is actually the only one that is invariant under the action of $\GL_4(\Z_p)$, see \cite[Theorem~1.1]{EL22}.
    \end{remark}
    
    \medskip
    
    For our computational experiments, we can determine the coefficients only up to finite precision. So in practice, we sample the variables $\xi_\alpha$ are considered to be random variables with uniform distribution on the set $\{ 0, 1, \ldots, p^{N + 1} - 1 \}$ for some precision $N$.

    \medskip
    
    We note that the Haar measure has already been used on the space of cubic surfaces over $\Q_p$ to compute the expected number of lines; see \cite{manssour2020probabilistic}. To the best of our knowledge, the other measures (to be defined next) are not recorded anywhere in the literature.

\subsection{The Blow-up Measure}

    As we saw in \cref{Sec:2}, one way to construct a cubic surface is to blow the plane $\P^2$ at six points in general position. If the six points are chosen randomly, we then get a random cubic surface. More precisely we define a new measure on the space of cubic surfaces as follows:
    
\begin{definition} (The blow-up measure)
    Let $\xi_0, \dots, \xi_6$ be independent uniformly distributed random variables in $\Z_p$ and let $F$ be the random polynomial $F = \xi_0 + \xi_1 X + \dots + \xi_6 X^6$. With probability $1$, the polynomial $F$ has degree $6$, its roots $\theta_1, \dots, \theta_6$ are simple and lie in the complement of the hyperplane arrangement \eqref{eq: arrangement}. This defines a measure on the space  of degree $6$ polynomials. The map $ F \longmapsto \SF $ determines a measure on the space of smooth cubic surfaces via pushforward which we call the \emph{blow-up measure}.
\end{definition}

    \medskip

    Again, in practice, sampling a polynomial $F$ is done by choosing a large integer $N$ and sampling the coefficients of $F$ from the set $\{ 0, 1, \ldots, p^{N + 1} - 1 \}$ independently with respect to the uniform distribution. The number of lines on the cubic surface $\SF$ is determined by factorizing $F$ over the field $K$ and applying Lemma~\ref{lem:NumberOfLines}.

    \medskip

    As mentioned in Remark~\ref{rem:DoubleSix}, not all cubic surfaces over $\Q_p$ arise from a blow-up of $6$ points in $\P^2$. Therefore, the blow-up measure only sees cubic surfaces with a Galois-invariant double-six.

\subsection{The Tropical (Generic) Measure}

    The following theorem suggests that the number of lines on a smooth cubic surface $\SF$ is tightly linked to the tropicalization of $F$:

    \begin{theorem}[Theorem 3.5 in \cite{OctanomialModel}] \label{thm:Octanomial_Theorem}
        Fix a prime $p \geq 5$ and a cubic surface $\SF$ over $\Q_p$, as in Definition~\ref{def: octanomial surface}. If $\SF$ is ``tropically smooth'', then the $27$ lines on $\SF$ have distinct tropicalizations in the tropical projective space $\mathbb{T}\P^3$. In that case, all $27$ lines on $\SF$ are defined over $\Qp$.
    \end{theorem}
    
    When it comes to the combinatorial type of the cubic equation $f$ (i.e., the regular subdivision it induces on the polytope $3 \Delta_3$), sampling from the Haar measure is not the best way of sampling polynomials with a diverse combinatorial types (or \emph{tropicalizations}). Under the Haar measure on $\Z_p$, it is quite rare to see elements with a large valuation. In order to remedy this shortcoming, we may sample from $\mathbb{Q}_p[ x_0 , x_1, x_2, x_3]_{( 3 )} $ by prescribing the valuation of each coefficient of $f$. More precisely,

    \begin{definition}[The tropical measure]
        Let $N$ be a positive integer and let $(\nu_{\alpha})_{|\alpha| = 3}$ be independent uniform random variables in $\{0, \dots, N\}$. We then obtain the random cubic equation
        \[
            f = \sum_{|\alpha| = 3} p^{\nu_\alpha} x^\alpha. 
        \]
        When $N$ is large enough, the cubic surface $\mathcal{S}(f)$ defined by $f$ is almost surely smooth. This defines a measure on the space of smooth cubic surfaces which we call the \emph{tropical measure}.
    \end{definition}

    It is quite natural to adapt this measure in the following way:

    \begin{definition}[The tropical generic measure]
        Let $N$ be a positive integer and let $(\nu_{\alpha})_{|\alpha| = 3}$ be independent uniform random variables in $\{0, \dots, N\}$ and $(c_{\alpha})_{|\alpha| = 3}$ be independent uniform random variables in $\Z_p^{\times}$ (with respect to the Haar measure on $\Z_p$). We then obtain the random cubic equation
        \[
            f = \sum_{|\alpha| = 3} p^{\nu_\alpha} u_\alpha  x^\alpha.
        \]
        The cubic surface $\mathcal{S}(f)$ is then smooth with probability $1$. The measure so defined on the space of smooth cubic surfaces is called the \emph{tropical generic measure}. 
    \end{definition}
    
\subsection{Experiments for \texorpdfstring{$p=7$}~}

    For each probability measure $\mu$ on the space of cubic surfaces, we can obtain a probability measure $\pi^{(\mu)}$ on the set $\{0,1,2,3,5,7,9,15,27\}$ of possible line counts, where $\pi^{(\mu)}_{i}$ records the probability under $\mu$ that a cubic surface contains $i$ lines. In general, given the measure $\mu$, it is quite hard to determine the distribution $(\pi^{(\mu)}_{i})$, even for the probability distributions defined in \cref{Sec:3}. Hence, we sample a large number of cubic surfaces under each measure and use Monte--Carlo estimation to get an idea of how this distribution looks like for the measures we defined above.

    \medskip
    
    We investigated each of the above measures experimentally by sampling  $10^5$ instances of smooth cubic surfaces and counting the corresponding number of lines. The latter is accomplished by using  Gr\"obner basis techniques, or rather Lemma~\ref{lem:NumberOfLines} for the blow-up measure. The resulting distributions on the number of lines are depicted in Table~\ref{tab:p-adicResults}. The code we used for sampling surfaces and counting lines can be found at (\ref{eq:link_for_code}).

    \begin{table}[H]
    \centering
    \begin{tabular}{|c|c|c|c|c|}
    \hline
    \# of $\Q_7$-Lines & Haar Mes. & Blow-up Mes. & Trop. Mes. & Trop. Gen. Mes. \\[0.6ex]
    \hline
     $0$  & 0.43995 &  0.19446 & 0.22230 & 0.25580 \\[0.6ex]
    \hline
    $1$  & 0.34534 & 0.12608 & 0.29004& 0.26891 \\[0.6ex]
    \hline
    $2$ & 0.08686 & 0.19604 & 0.02106& 0.02083 \\[0.6ex]
    \hline
    $3$ & 0.08564 &  0.21602 & 0.29145& 0.26952 \\[0.6ex]
    \hline
    $5$  &0.03169 &0.12778 & 0.04620 & 0.04988 \\[0.6ex]
    \hline
    $7$  & 0.00467 & 0.07337 &0.06708 & 0.06911\\[0.6ex]
    \hline
    $9$  & 0.00401& 0.05099 & 0.02868& 0.03443\\[0.6ex]
    \hline
    $15$  & 0.00045 & 0.01500 & 0.02582& 0.02666\\[0.6ex]
    \hline
    $27$  & 0.00001 & 0.00026 &0.00487 & 0.00406\\[0.6ex]
    \hline
    \textbf{Average} & 1.01023 & 3.00964 & 2.68398 & 2.67200\\[0.6ex]
    \hline
    \end{tabular}
    \caption{The distribution of the number of lines for different probabilistic measures over the $7$-adic numbers.}
    \label{tab:p-adicResults}
    \end{table}

    It is important to note the following:
    
    \begin{enumerate}
    
        \item Our result only estimates\footnote{Using a Monte-Carlo method.} the distribution $\pi^{(\mu)} = \left(\pi^{(\mu)}_{k}\right)$ for any chosen measure $\mu$. This approximation relies on the law of large numbers, so our results are random but converge to the correct distribution as the sample size gets larger and larger.

        \medskip

        \item The dimension of the space of cubic equations is $20$, which would make the computation of integrals\footnote{The probability that a random surface has a certain number of lines is an integral.} quite expensive. One of the advantages of the estimation method we used is that it is not affected by the dimension of the space we sample from.

        \medskip

        \item When dealing with computations over a $p$-adic field, one has to be careful with precision. Our computation were conducted with an absolute $p$-adic precision of 300, while a random $p$-adic number in $\Z_p$ was sampled from the uniform distribution on the set $\{0,\dots,p^{8}-1\}$. All these choices were made so that the computation runs in a reasonable time all the while keeping the results significant. 

        \medskip

        \item Finally, it should also be mentioned that our implementations in \eqref{eq:link_for_code} can be optimized a great deal from a performance point of view. So it is very much possible to go beyond the sample sizes we have used.

    \end{enumerate}

\subsection{Interpretation of the Results}

    We observe from Table~\ref{tab:p-adicResults} that under the Haar measure the probabilities of the line counts decrease. Actually most surfaces have $0$ or $1$ $\Q_p$-lines so under this measure, it is quite rare to find a cubic surfaces with a high count of $p$--adic lines, in particular $27$-lines (around $\simeq 10^{-5}$ probability). This explains that the average number of $p$-adic lines under the Haar measure is almost $1$.
    In fact, the following theorem quantifies the expected number which is $\simeq 1.01749$ for $p = 7$.
    
     \begin{theorem}\emph{(\cite[Theorem~3]{manssour2020probabilistic})}
     \label{thm: manssour 2020} 
        The expected number of $p$-adic lines on a random uniform $p$-adic cubic surface in $\P^3$ is $\frac{(p^3-1)(p^2+1)}{p^5-1}$.
    \end{theorem}

	We observe that under both tropical measures, the chance of seeing bigger line counts is significantly higher compared to the Haar measure. Notice also that both measures yield more or less the same distribution of line counts. This is probably due to the fact that the number of lines depends heavily on the tropicalization of the cubic equation $f$, and hence also on the induced triangulation on the newton polytope $3\Delta_3$. This is in light of \cite[Conjecture~4.1]{OctanomialModel}.
    
    As far as the blow-up measure is concerned, we see again bigger number of lines with higher probability compared to any of the other measures. This can be explained by the fact that this measure only sees the cubic surfaces with a $\Q_p$-rational double-six, which generally have more lines than those that do not. This is in line with the results from \cite[Theorem~1]{Manjul2021} where the probability $\rho_n(r)$ that a degree $n$ polynomial with random independent and uniform coefficients in $\Z_p$ has exactly $r$ roots is given recursively. The particular case of interest to us is the case $n = 6$ and $p = 7$, where we get
    \begin{align*}
        \rho_6(0) + \rho_6(1) &= \frac{7280010099060058135701356421229303451929}{9915124900168002703437229470076926702000} \simeq 0.7343, \\
        \rho_6(2) &= \frac{132142086852025648305980401844338671377}{661008326677866846895815298005128446800} \simeq 0.19991, \\
        \rho_6(3) &= \frac{6274737460539590192834937928502919283}{123939061252100033792965368375961583775} \simeq 0.05063, \\
        \rho_6(4) &= \frac{3279805090966404942802616745034685803}{220336108892622282298605099335042815600} \simeq 0.01489,\\
        \rho_6(5) &= 0,\\
        \rho_6(6) &= \frac{379252487878267254806025930638752849}{1101680544463111411493025496675214078000} \simeq 0.00035.
    \end{align*}
    \cref{tab:compareManjul} compares the exact results from  \cite{Manjul2021} to the estimates obtained in our experiments.

    	\begin{center}
    		\begin{table}[H]
	    	\begin{tabular}{|c|c|c|c|c|c|}	   
				\hline
				\# of roots of $F$  &  $ 0,1$  	 &  $ 2$  &  $3$     &    $4$   &    $6$        \\[0.6ex]
	    		\hline
				\# of lines on $\SF$  &   $0,1,2,3$    & $5,7$    &  $9$     &    $15$   &    $27$        \\[0.6ex]
				\hline	
				 Results from \cite{Manjul2021}              & 0.73423     &  0.19991      &   0.05063     &  0.01489    & 0.00035     \\[0.6ex]
				\hline
				Estimates from \cref{tab:p-adicResults}&  0.73260            &     0.20115       &   0.05099    &    0.01500    &    0.00026     \\[0.6ex]
				\hline
			\end{tabular}
			\caption{Comparison of theoretical and experimental results for the blow-up measure.}
			\label{tab:compareManjul}		
    		\end{table}
    \end{center}
 
    \section{Heuristics for the Real Numbers}\label{Sec:4}
    
    Let $V = \R[x_0,x_1,x_2,x_3]_{(3)}$ be the real vector space of homogeneous degree three polynomials in $x = (x_0,x_1,x_2,x_3)$. We endow $V$ with the inner product $\langle \cdot, \cdot \rangle$ defined as
    \[
        \langle f_1,f_2 \rangle \coloneqq   \frac{1}{4\pi^2} \int_{\R^{4}} f_1(x) f_2(x) e^{- \frac{\norm{x}^2}{2}} dx.
    \]
    The space $V$ is endowed with the action of the orthogonal group $O(4,\R)$ by change of variables, i.e.,
    \[
        (g\cdot f)(x) = f(g^{-1}x) \quad \text{for any } f \in V, \ g \in O(4,\R).
    \]
    This makes $V$ a linear representation of $O(4,\R)$ which is moreover unitary, i.e., the inner product $\langle \cdot , \cdot \rangle$ is stable under $O(4,\R)$:
    \[
        \langle g \cdot f_1 , \  g\cdot f_2 \rangle = \langle f_1 , f_2 \rangle \quad \text{ for } f_1,f_2 \in V \text{ and } g \in O(4,\R).
    \]
    However, $V$ is not irreducible, so $\langle \cdot, \cdot \rangle$ is the not the unique $O(4,\R)$-invariant inner product on $V$. The representation $V$ splits into two irreducible orthogonal sub-representations as follows:
    \begin{equation}\label{eq: Harmonica}
        V = \mathcal{H}_3 \oplus \norm{x}^2 \cdot \mathcal{H}_1
    \end{equation}
    where $\mathcal{H}_3$ is the subspace of homogeneous harmonic polynomials of degree $3$ and $\mathcal{H}_1$ is the space of homogeneous degree $1$ polynomials. Let $p_1$ and $p_2$ be the orthogonal projections of $V$ onto $\mathcal{H}_3$ and $\norm{x}^2\mathcal{H}_1$, respectively.
    \begin{definition}
        \label{def: invariant inner products}

            Let $\lambda$ and $\mu$ be positive real numbers. We define the inner product $\langle \cdot, \cdot \rangle_{\lambda, \mu} $ on $V$ as follows:
    \[
        \langle f_1 , f_2 \rangle_{\lambda, \mu}  = \frac{1}{\lambda^2} \ \langle p_1(f_1), p_1(f_2)\rangle \  + \  \frac{1}{\mu^2}  \ \langle p_2(f_1), p_2(f_2)\rangle.
    \]
    We also define the centered Gaussian measure $(\P_{\lambda, \mu})_{\lambda, \mu > 0}$ on $V$ associated to $\langle \cdot, \cdot \rangle_{\lambda, \mu}$ i.e. the measure on $V$ whose density $\phi_{\lambda, \mu}$ is proportional to
    \[
      \phi_{\lambda, \mu} (f) \propto \exp\left(\frac{- \langle f , f \rangle_{\lambda, \mu}}{2}\right), \quad  \text{for } f \in V.
    \]
    
    \end{definition}

    The two parameter family in \cref{def: invariant inner products} parametrizes all $O(4,\R)$-invariant inner products on $V$ and hence all non-degenerate Gaussian measures on $V$ that are $O(4,\R)$-invariant are of the form $\P_{\lambda, \mu}$ (see \cite[Section~4]{Kost90} for more details).
    
    \medskip
    
\begin{figure}[H]
        \centering
        \includegraphics[scale=0.7]{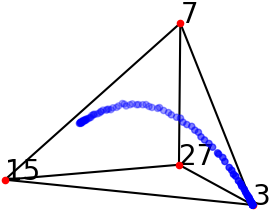}
        \caption{The curve $(\pi^{(\lambda)})_{\lambda \in (0,1)}$ inside the $3$-dimensional probability simplex.}
        \label{fig:curveReals}
    \end{figure}

    To sample a random cubic $f \in V$ with respect to the measure $\P_{\lambda, \mu}$, we use the orthonormal basis $\{H_{3,i}\}_{1 \le i \le 16} $ and $ \{H_{1,j}\}_{1 \le j \le 4}$ respectively of $\mathcal{H}_3$ and $\norm{x}^2 \mathcal{H}_1$ given in \cite[Table~1]{realLines}. More precisely, let $(\xi_{3,i})_{1 \leq i \leq 16}$ and $(\xi_{1,j})_{1 \leq j \leq 4}$ be two independent sequences of independent identically distributed standard Gaussian random variables and let us define the random cubic $f$ as follows:
    
    \[
        f =  \lambda \sum_{1 \le i \le 16} \xi_{3, i} H_{3, i}  \ + \ \mu \sum_{1 \le j \le 4} \xi_{1, j}  H_{1, j}.
    \]
    Since multiplying with a scalar does not change the zero set of $f$, we may assume that  $\lambda + \mu = 1$. Then we may focus on the $1$-parameter family of measures
    \[
        \P_\lambda \coloneqq \P_{\lambda, 1 - \lambda}, \quad \lambda \in (0,1).
    \]
    For each $\lambda \in (0,1)$, we have sampled $10^5$ cubics under $\P_{\lambda}$ and estimated the probability distribution $\pi^{(\lambda)} = \big(\pi^{(\lambda)}_3, \pi^{(\lambda)}_7, \pi^{(\lambda)}_{15}, \pi^{(\lambda)}_{27} \big)$ where
    \[
    		\pi^{(\lambda)}_k = \P_{\lambda} \left(  \#\{ \R\text{-rational lines in } \mathcal{S}(f)\} = k \right) \quad \text{ for } k = 3,7,15,27.
    \]
    The result is a curve 
    \[
        (0,1) \to \Delta_3, \quad \lambda \mapsto \pi^{(\lambda)} = \big(\pi^{(\lambda)}_3, \pi^{(\lambda)}_7, \pi^{(\lambda)}_{15}, \pi^{(\lambda)}_{27} \big),
    \]
    in the $3$-dimensional $\Delta_3$ simplex:
    \[
        \Delta_3 \coloneqq \big\{ (p_1, p_2, p_3, p_4) \in \R_{\geq 0}^4 \colon  p_1 + p_2 + p_3 + p_4 = 1\big\}.
    \]
    The $3$-dimensional simplex $\Delta_3$ is contained in the affine hyperplane in $\R^4$ given by $\{x_1 + x_2 + x_3 + x_4 = 1\}$. We pick an orthogonal basis $(e_1, e_2, e_3)$ of the hyperplane $\{x_1 + x_2 + x_3 + x_4 = 0\}$ and write expand the curve $\lambda \mapsto \pi^{(\lambda)}$ as follows

    \[
        \pi^{(\lambda)} = \frac{1}{4} (1,1,1,1) \  +  \  y_{1}(\lambda) e_1 + y_{2}(\lambda) e_2 + y_{3}(\lambda) e_3.
    \]
   \cref{fig:curveReals} depicts the curve $\lambda \mapsto (y_{1}(\lambda), y_{2}(\lambda), y_{3}(\lambda))$.
    
	A particular case of interest is the distribution $\pi^{(1/3)} $ obtained when we sample from the Kostlan distribution $\P_{1/3}$. To estimate $\pi^{(1/3)}$, we sampled $10^6$ cubics $f \in V$ under $\P_{1/3}$ and counted the number of lines on $\mathcal{S}(f)$. Table~\ref{table:KostlanMes} summarizes the estimate we obtained for $\pi^{(1/3)}$.
	
	\begin{center}
		\begin{table}[H]
				\begin{tabular}{|c|c|c|c|c|}
					\hline
				  $\pi^{(1/3)}_{3}$ &   $\pi^{(1/3)}_{7}$ &   $\pi^{(1/3)}_{15}$ &   $\pi^{(1/3)}_{27}$ & Average \# of lines  \\[0.6ex]
					\hline
			         0.570252         &    0.338973            &     0.089406          &  0.001369  &  5.46162  \\[0.6ex]
					\hline	 
				\end{tabular}
			\caption{The distribution of number of real lines under the Kostlan measure.}
			\label{table:KostlanMes}
		\end{table}
	\end{center}
    
    We note that our experiments corroborate the following results of \cite{BasuLeri} and \cite{realLines}:
    
    \begin{theorem}\emph{(\cite[Theorem~5]{BasuLeri})}
        The average number of real lines on a random cubic surface in $\R\P^3$ under the Kostlan distribution $\P_{\frac{1}{3}}$ is $6\sqrt{2} - 3 \simeq 5.48528$.
    \end{theorem}

    \begin{theorem}\emph{(\cite[Theorem~1]{realLines})}
        The expected number of real lines on a random real cubic surfaces under the measure $\P_{\lambda}$ is 
        \[
            E_{\lambda} = \frac{9(8 \lambda^2 + (1-\lambda)^2)}{ 2\lambda^2 + (1-\lambda)^2} \left( \frac{ 2\lambda^2 }{ 8 \lambda^2 + (1-\lambda)^2 } - \frac{1}{3}  +  \frac{2}{3}\sqrt{  \frac{8 \lambda^2 + (1-\lambda)^2}{ 20 \lambda^2 + (1-\lambda)^2} } \right).
        \]
    \end{theorem}
    Our computations were carried out using the Julia package HomotopyContinuation.jl \cite{HomotopyContinuation.jl}.  Code and data are available at \eqref{eq:link_for_code}.

    \subsection{Interpretation of the Results}
        
    	Notice that, as $\lambda \to 0$, the distribution $\pi^{(\lambda)}$ converges to the vertex $ \pi^{(0)} = (1,0,0,0)$. So when $\lambda$ is small, we can only hope to see surfaces with $3$ real lines under $\P_{\lambda}$; see \cite[Proposition~2]{realLines}. As $\lambda \to 1$ we see higher number of lines with bigger probabilities. The curve in \cref{fig:curveReals} depicts of $\pi^{(\lambda)}$ as $\lambda$ ranges in $(0,1)$.
        We observe that the limit distribution as $\lambda \to 1$ (a random cubic surface under the measure $\P_{1}$ supported on $\mathcal{H}_3$ is smooth with probability $1$) lies in the interior of the $3$-dimensional probability simplex but is not easy to determine explicitly.
        
        This poses the following question:

    \begin{question}
        Is the curve $(\pi^{(\lambda)})_{\lambda \in (0,1)}$ algebraic, and if yes, what are the equations defining it? 
    \end{question}

    \section{Galois Groups} \label{Sec:5}

    Let $S$ be a smooth cubic surface defined over a perfect field $K$, and let $L$ be the field of definition of the $27$ lines on $S$. Then the extension $L/K$ is Galois, and the Galois group $\Gal(L/K)$ is a subgroup of $W(E_6)$, the Weyl group of order $51840 = 2^7\cdot 3^4\cdot 5$. A natural question to ask is which subgroups of $W(E_6)$ can be realized in this way. The answer depends, of course, on the base field $K$. Elsenhans--Jahnel proved the following:
    \begin{theorem}\emph{(\cite[Theorem~0.1]{ElsenhansJahnel2015})}
    \label{thm:Elsenhans--Jahnel}
    All subgroups of $W(E_6)$ can be realized when $K=\Q$.
    \end{theorem}

    Let us now consider the case where $K$ is a non-archimedean local field of characteristic $0$; such a field is isomorphic to a finite extension of $\Q_p$ for some prime $p$. In this case, there are certain constraints on the possible Galois groups that can arise. Let us make this more precise. Let $L$ be a finite Galois extension of $K$ with Galois group $G$, and consider (the first parts of) the ramification group filtration on $G$:
    \[G \supseteq G_0 \supseteq G_1.\]
    Here, $G_0$ and $G_1$ are the inertia and ramification groups, respectively. It is well-known that
    \begin{itemize}
        \item the group $G_1$ is a $p$-group,
        \item the quotient $G/G_0$ must be cyclic, and
        \item the quotient $G_0/G_1$ must be cyclic of order coprime to $p$
    \end{itemize}
    (see, for example, \cite[Chapter~IV]{SerreLocalFields}). In addition to these constraints, we also have the following:
    \begin{proposition}
        The Galois group $G$ is solvable.
    \end{proposition}
    \begin{proof}
        This well-known result follows from the fact that every $p$-group is nilpotent.
    \end{proof}
    
    Because of this proposition, the statement in Theorem~\ref{thm:Elsenhans--Jahnel} clearly can not be true for $K$ since $W(E_6)$ has non-solvable subgroups. Indeed, there are precisely $19$ non-solvable subgroups of $W(E_6)$ up to conjugation, and they can be computed using the following Magma code:
    
    \medskip

    \begin{verbatim}
        R_E6 := RootDatum("E6");
        Cox_E6 := CoxeterGroup(R_E6);
        WE6 := StandardActionGroup(Cox_E6);
        print NonsolvableSubgroups(WE6);
    \end{verbatim}
    These observations pose the following question:
    \begin{question}
    \label{question: which can be realized}
    Which subgroups of $W(E_6)$ can arise as Galois groups of lines on smooth cubic surfaces  over $\Q_p$? How does this list depend on $p$?
    \end{question}

     The results of our experiments suggest that Galois groups should be quite small and they are usually abelian. Tables~\ref{table:abeliangroups_p=5}-\ref{table:nonabeliangroups_p=5} (resp. Tables~\ref{table:abeliangroups_p=7}-\ref{table:nonabeliangroups_p=7}) summarize the Galois groups\footnote{In the tables, we use the notation provided by Magma. In particular, $F_q$ is the Frobenius group $\F_q \rtimes \F_q^{\times}$, and $OD_{2^k}$ is the ``other-dihedral'' group $C_{2^{k-1}} \rtimes C_2$ with $C_2$ acting as $2^{k-2} + 1$.} obtained for a sample of $25000$ surfaces over $\Q_5$ (resp. $\Q_7$) sampled from the Haar measure, and the number of times they occurred. Our code can be found at (\ref{eq:link_for_code}).
	\begin{center}
		\begin{table}[H]
			\begin{tabular}{|c|c|c|c|c|c|c|c|c|c|c|c|c|c|c|c|c|c|c|c|c}
				\hline
				$C_2$ &   $C_3$ &   $C_4$ &   $C_{5}$  &   $C_{6}$ &  $C_{8} $  & $C_9$ & $C_{10}$ & $C_{12}$ &  $C_{2}^2$ &   $C_{4}^2$ &   $C_{2}\times C_{4}$  &   $C_{2}\times C_{6}$   \\[0.6ex]
				\hline
				242           &    246             &     2117           &  2042   &  5154 & 2928 & 2126   & 	2647           &    3613             &     475           &  30   &  1193 & 1294   \\[0.6ex]
				\hline	  
			\end{tabular}
			\caption{The $24107$ abelian groups that appeared in our sample for $p=5$.}
			\label{table:abeliangroups_p=5}
		\end{table}
		\begin{table}[H]
			\begin{tabular}{|c|c|c|c|c|c|c|c|c|c|c|c|c|c|c|c|c|c|c|c}
				\hline
				$D_5$ & $D_6$ & $D_{10}$ & $F_5$ & $S_3$ & $OD_{16}$ & $C_{2}\times F_{5}$ & $C_{3}\times S_{3}$   & $C_{4}\times S_{3}$ & $C_{6}\times S_{3}$ & $C_3\rtimes C_4$ & $C_3\rtimes C_8$    \\[0.6ex]
				\hline
				2 & 98 & 2 & 13 & 204 & 7 & 21 & 323 & 4 & 47 & 170 & 2  \\[0.6ex]
				\hline	  
			\end{tabular}
			\caption{The $893$ non-abelian groups that appeared in our sample for $p=5$.}
			\label{table:nonabeliangroups_p=5}
		\end{table}
		\begin{table}[H]
			\begin{tabular}{|c|c|c|c|c|c|c|c|c|c|c|c|c|c|c|c|c|c|c|c|c|c}
				\hline
				$C_2$ &   $C_3$ &   $C_4$ &   $C_{5}$  &   $C_{6}$ &  $C_{8} $  & $C_9$ & $C_{10}$ & $C_{12}$ &  $C_{2}^2$ &   $C_{3}^2$ &   $C_{2}\times C_{4}$  &   $C_{2}\times C_{6}$ & $C_{3}\times C_{6}$  \\[0.6ex]
				\hline
				284           &    286             &     2263           &  2149   &  5513 & 3001 & 2303   & 	2835           &    3931             &     333           &  42   &  787 & 1043 & 152 \\[0.6ex]
				\hline	  
			\end{tabular}
			\caption{The $24922$ abelian groups that appeared in our sample for $p=7$.}
			\label{table:abeliangroups_p=7}
		\end{table}
		\begin{table}[H]
			\begin{tabular}{|c|c|c|c|}
				\hline
				$D_4$ & $F_5$ & $C_2\times F_5$ & $C_4\rtimes C_4$    \\[0.6ex]
				\hline
				50 & 7 & 3 & 18   \\[0.6ex]
				\hline	  
			\end{tabular}
			\caption{The $78$ non-abelian groups that appeared in our sample for $p=7$.}
			\label{table:nonabeliangroups_p=7}
		\end{table}
	\end{center}
	
    \begin{remark}
            Notice that the Galois groups that appeared for $p = 5$ are more complicated. This is expected since the prime $5$ divides the order of the Weyl group $W(E_6)$. We have also tried to make the same computation for the other primes with the same property, namely $2$ and $3$. Determining Galois groups using the Magma's \texttt{GaloisGroup()} function, however, was significantly slower in these cases. This means that more interesting groups show up for the primes $2$ and $3$.
    \end{remark}

    \bigskip

	We see that while the generic Galois group for surfaces over $\Q$ is $W(E_6)$, over the $p$-adics there is no ``generic'' group but rather a list of ``small'' groups that can occur with positive probability. Question \ref{question: which can be realized} would be interesting to answer in future work. We conclude this section with the following remark.
 \begin{remark}
     Elsenhans--Jahnel also showed that some Galois groups are possible over any field of odd characteristic as long as a field extention with that group exists; see \cite{ElsenhansJahnel2019a,ElsenhansJahnel2019b}. This might give some information on what happens over $\Q_p$ by lifting surfaces over $\F_p$ to $\Z_p$.
 \end{remark}

    \bibliographystyle{alpha}
    \bibliography{refs}

\end{document}